\documentclass{amsart}

\begin{document}

\title{Aristotle's Relations:
An Interpretation in Combinatory Logic}     
\author{Erwin Engeler}
\address{ETH Z\"urich}   
\date{\today} 
\maketitle

\textsc{Abstract}
\\
The usual modelling of the syllogisms of the Organon by a calculus of classes does not include
relations. Aristotle may however have envisioned them in the first two books as the
category of relatives, where he allowed them to compose with themselves. Composition
is the main operation in combinatory logic, which therefore offers itself to logicians for a
new kind of modelling. The resulting calculus includes also composition of predicates by
the logical connectives.

\vspace{10 mm}

\textbf{Introduction}

\vspace {3 mm}

Relations turn up at birthdays,%
\footnote{
The first part of his little essay was written as a gift for my 90th birthday, (see dedication).
}
congratulating. Even logicians have them; of the first one,
Aristotle, we even know some of their names.
But it is a question among historians of axiomatic geometry whether he had the other
kind of relations, the ones that modern logicians are concerned with. Indeed, many hold
that the Stagirite did not have this concept, and that Greek mathematics, in particular Euclid
shows this: the relation of betweenness does not enter his axioms for geometry. The
tradition of Aristotelian Logic is often blamed for this serious lacuna. In fact Euclidean
axiomatics reached completion only in the 19th century, [4]. Missing the concept of relation
is perceived as making it unfit as an adequate logic for developing formal axiomatic
mathematics; that had to wait till the Boole-Peano-Russell disruption which eclipsed traditional
logic. This, I think, is overstating the case. The Organon itself gives a picture of
Aristotle's understanding of judgements other than those that are formulated by the syllogistics
and included relations and functionals.
My argument comes in two parts, motivation and formal development.
The first part experiments with concepts of definitions on an example that Aristotle himself
could have handled. These are discussed in terms, motives, that I discern in the
Organon, in particular the composition of predicates called \textquotedblleft relatives", the use of logical
connectives and forms of recursion. Modalities, also an important ingredient of the
Organon are not included here, which while feasible are irrelevant to the present theme.
In the second part, the compositional aspect of predicates is brought forward and made
into the basis for interpreting the {\em Organon}. This results in the establishment of a logical
calculus $\mathcal{E}^\Lambda$ of judgements. This provides a model of syllogisms which include relations.
The conclusion is, that Aristotle had the means to treat relations but chose not to do so for
his syllogistics.

 \newpage
\textbf{ 1 \quad Aristotle'  Relatives}

\vspace { 3 mm}

If you'll bear with me, let us see how Aristotle would, and could speak of his relatives, formally, and within the framework of his toolkit, the \emph{Organon }--- in my naive and quite ahistoric reading, using an \emph{ad hoc } formalism that we shall later turn into a formal calculus. 

\vspace{ 3 mm}

\textsc {1.1 Predication}
\\
The main grammatical operation is applying a predicate to a subject: $[\textit{red}] \cdot [\textit{blood}] $ is a statement which predicates that blood is red. All kinds of thought objects are admitted as predicates; Aristotle divides them into \textquotedblleft categories", distinguishing for example between predicates about quantity (\textquotedblleft big"), and quality (\textquotedblleft red") and relatedness, (the category of \textquotedblleft relatives").
\\
Looking at example of relatives, let [\textit{mother}] predicate of a subject that it is enjoying motherhood. Thus $[\textit{mother}] \cdot  [\textit{Phastia}]$ states that Phastia is a mother. Nothing prevents us from using this thought object as a predicate: 
   \\
  $ ( [\textit{mother}] \cdot  [\textit{Phastia}] ) \cdot [\textit{Aristotle}]$ 
\\
tells us that Phastia is the mother of Aristotle. This turn is what Aristotle (Categories, Chapt.7) really meant with the category of relatives; he in fact called this category \textquotedblleft things pointing towards something". The above predicate $[\textit{mother}]$ of motherhood  predicates of a subject $a$ that the subject $b$ is her child, \textquotedblleft pointing $a$ to $b$". Thus, a relative predicate in fact introduces a binary relation by composition of predicates. The compositional nature of relatives is an early showing of Currying, a device that became one of the central aspects of combinatory logic.
\\
For the moment, we don't concern ourselves with the question as to what category some predicate or subject might belong, all are treated as relatives; one of the uses that Aristotle gets out of such prescriptions is to avoid predicating nonsense by disqualifying predications between certain categories.
\\
Using the \textquotedblleft relative" predications of motherhood and fatherhood, of marriage, and of being male or female applied to family members, we can easily envision the genealogy of Aristotle as a list of such statements. He himself would use many more predicates to talk about his relations: he would use \textquotedblleft son", \textquotedblleft sister", \textquotedblleft grandfather", \textquotedblleft sister-in-law", or even \textquotedblleft male descendent" etc. as predicative concepts. Let us see how that could fit in.

\vspace{3 mm}

\textsc{1.2 Explicit Definitions}

\vspace{ 3 mm}

One way of introducing a new predicative concept is explicitly, as composite predicates. Using variables $a,b,c, \dots$ for subjects, the definition of $b$ being a child of $a$ is simply
\\
$ ([\textit{child}] \cdot b ) \cdot a    = ([\textit{mother}] \cdot a ) \cdot b $. 
\\
For some relationships we need logical connectives such as \textquotedblleft and" and \textquotedblleft or". These are denoted by $\wedge$ and $\vee$, and used on predicates $P$ and $Q$  to obtain $P \wedge Q$ and $P \vee Q$. Thus, the predicate of being a \textquotedblleft son" is
\\
 $([\textit{son}] \cdot b ) \cdot a    = ([\textit{mother}] \cdot a ) \cdot b \wedge [\textit{male}]b$,
 \\
Similarly, constituting the predicative concept of a \textquotedblleft family": For $a,b,c,d$  to form a core family, predicated by the predicate $[family]$ on some individuals $a, b, c, d$ set
\\
$[\textit{family}]abcd  = [\textit{mother}]ac  \wedge  [\textit{mother}]ad \wedge   [\textit{father}]bc \wedge   [\textit{father}]bd \wedge [\textit{female}]a \wedge [\textit{male}]b$.
\\
\emph{Notation}: we have dropped the center-dot that denotes application and adhere to the convention that sequences of applications are to be understood as parenthesised to the left: $uvw$ is read as $(u \cdot v) \cdot w$, etc.
\\
The idea of \textquotedblleft pointing to something" has just been applied again: if $Pa$ is a relative predication then $Pa$ and $Pab$ may be too. $P$ in the context $Pabc$, for example, would introduce a ternary relation. This leads to a more general notion of explicit definitions: 
\\
Formally, an \emph{explicit definition} concerns an expression $\varphi(x_1, x_2, \dots x_n)$ built up from variables, predicates (introduced earlier as definienda) and the logical connectives. It defines a new predicate $P$ by the defining equation
 \\
 $Px_{1} x_{2 } \dots x_{n} = \varphi(x_1, x_2, \dots, x_n)$.
\\
This is the \emph{Principle of Comprehension}. It comprehends the connective structure of $\varphi$ into a single predicate, an idea that goes back to Sch\"onfinkel and is a basic concept for Curry's combinatory logic.

\vspace{3 mm}

\textsc{1.3 Implicit Definitions}

\vspace { 3 mm}

As ethnologists tell us, all, (even \textquotedblleft primitive") cultures allow definitions of familial relatedness, some quite elaborate. The simplest ones are of the explicit kind as above, but this is far from sufficient. Consider the notion of being (maternal) siblings. The desired predicate $[sibling]cd$ hides a mother somewhere in its definition. Aristotle resolves this by introducing the construction \textquotedblleft some $P$".
\\
We denote the construct \textquotedblleft some $P$" by $\varepsilon_{x}(P x)$. It implies a sort of existential referent. The assertion \textquotedblleft some $P$ are $Q$" would then transcribe to  $Q \cdot (\varepsilon_{x}(P x))$. This device was introduced by Hilbert to be used in his foundational program as a tool for the formal elimination of quantifiers. With it we can define the predicate of being a sibling:
\\
 $[\textit{sibling}]yz = [\textit{mother}] (\varepsilon_{x} ( [\textit{mother}] x y )) z $,
\\ 
stating that $y$'s mother $x$ is also mother of $z$. 

\vspace{3 mm}

\textsc{1.4 Recursion}

\vspace{ 3 mm}

We now can have uncles, cousins of  all kinds, marriages  between them, etc., enough to tell the story of Aristotle's relations, and formulate things like \textquotedblleft Aristotle is the father-in-law of his niece". However, the concept  \textquotedblleft male descendent" which was very important at the time is still open. This is a typical case for using a definition by recursion: 
\\
$[\textit{mdesc}]xy = ([\textit{male}]y \wedge [\textit{father}]xy ) \vee (\varepsilon_{z} ([\textit{father}]zy \wedge [\textit{mdesc}]xz ))$
\\
stating that a male descendent $y$ of $x$ either is a son of $x$ or there is some male descendent who is his father. 
\\
Are such definitions admissible? What it amounts to, is that it allows the use of the $\varepsilon$-operator also for predicates as variables: The definiendum $[\textit{mdesc}]$ is represented by a variable $U$. 
\\
$[mdesc]xy = \varepsilon_{U}(([male]y \wedge [father]xy ) \vee (\varepsilon_{z} ([father]zy \wedge Uxz )))$
\vspace{3 mm}
\\
You may have noticed  that the $\varepsilon $ operator allows to understand the Aristotelean \textquotedblleft some $P$ are $Q$" as $Q \cdot (\varepsilon_{x} Px )$ in the present \emph{ad hoc} formalism. We have not used two other basic ingredients of Aristotle's formalism; there was no need of negation and of the companion to \textquotedblleft some $x$ are $y$", namely \textquotedblleft all $x$ are $y$". This will be be task of a later chapter.
\vspace{ 10 mm}

\textbf {2 \quad Discussion}

\vspace {3 mm}

\textsc{ 2.1 On the \emph{Organon}}

\vspace { 3 mm}

For all I know, Aristotle would have accepted each definiens in the above definitions as a statement in his sense. But he would have hesitated to call them \textquotedblleft categorical" statements. The distinction arises in the course of his development of the \emph{Organon}\footnote{Consulted translations: Owen [7] and Smith [10].}.
\\
The \emph{Organon}, the way I read it, has the character of a manual, a texbook that instructs the reader in the art of preparing a conclusive argument using well-formed and immediately understandable statements. The \emph{Organon} consists of three parts, \textquotedblleft The Categories", \textquotedblleft On Interpretation" and \textquotedblleft Prior Analytics" (in two books).\footnote{These are cited as (Cat), (Int), (PA A), (PA B) with only the chapters indicated, e.g. (Cat 8) is Chapter 8 in Categories.}
\\
The first two introduce the notion of well-formed-ness by a discourse of examples and grammatical distinctions. The question is how to combine predicates while evading non-sensical, ambiguous or misunderstandable compositions as far as possible:  How do you compose $[\textit{good}]$, $[\textit{man}]$ and $[\textit{shoemaker}]$ (Int  11) ? 
\\
You may take conjunction and application to form
\\
$([\textit{good}] \wedge [\textit{shoemaker}])[\textit{man}]$ or $([\textit{good}][\textit{shoemaker}])[\textit{man}]$ or $[\textit{shoemaker}]([\textit{good}] [\textit{man}])$ or $([\textit{good}] [\textit{man}])[\textit{shoemaker}]$,
\\
each with a valid and different meaning,.
\\
Disjunction can enter with predicates of the category of qualities such as colour; $([black] \vee [white])[man]$, (Cat 8), or
\\
$([\textit{green}] \vee [\textit{blue}]) [\textit{turquoise}]$, or $[\textit{green}] ([\textit{blue}] [\textit{turquoise}])$ or $([\textit{blue}] [\textit{green}]) [\textit{turquoise}]$,
\\
as you may judge the quality of the stone.
\\
Negation enters in two pairs of \emph{opposites } (Cat 10, Int 6-10):
\\
\textquotedblleft some $P$ are $Q$" as against \textquotedblleft no $P $ is $Q$", and
\\
\textquotedblleft all $P$ are $Q$" as against \textquotedblleft some $P$ are not $Q$".
\\
These are the four kinds of statements that Aristotle calls \emph{categorical statements}.

\vspace {3 mm}

In Prior Analytics, Aristotle gets down to the business of constructing conclusive sequences of arguments. Already required for categorical statements, each argument should be terse and immediately understandable. He chose to restrict the form of statements to the four categorical ones above, which for convenience we denote in the form $\varphi (x, y), \psi(x, y)$, etc.
\\
And then he shows \emph{unguem leonis}\footnote { \textquotedblleft Writing with the the claws of a lion."}: he develops a formal system of logic based on logical arguments, called \emph{syllogism}, of the form
\\
$\varphi(P, Q ), \psi(Q, R) \vdash \chi(P, R)$,
\\
expressing that from the categorical statements $\varphi(P, Q )$ and $ \psi(Q, R)$ you may infer $\chi(P, R)$, (PA  A 1 - 14). Since there are four of these, there is a plethora of such deductive patterns. Aristotle proceeds to eliminate all but fourteen of them by showing, using counterexamples, which of them preserve the truth of the statements.
\\
The high point of Prior Analytics is the proof of a metatheorem: all fourteen can be deduced from just two of them, (PA A 23).
\vspace{3 mm}

\textsc{2.2 Critiques of the \emph{Organon}}

\vspace { 3 mm}

By this metatheorem, Aristotle establishes a sort of completeness for his syllogistic proof system. This is different from the present notion of completeness in mathematical logic which involves models. 

\vspace { 3 mm}

The story of \emph{models} for Aristotle's logic can be traced from Boole's Laws of Thoughts, to \L ukasiewicz [6], Shepherdson  [9] and Corcoran [1] in the 20th century.
The view developed that this logic dealt with unary truth-functions which could be understood as classes. To accomodate the interpretation of categorical sentences, the classes had to be non-empty. Corcoran, for example, interpreted these as:
\\
\textquotedblleft all $ A $ are $ B$"  is: $A \subseteq B$, 
\\
\textquotedblleft no $A$ is $B$" is:  $A \cap B = \emptyset $,
\\
\textquotedblleft some $A$ are $B$" is: $A \cap B \neq \emptyset$,
\\
\textquotedblleft some $A$ are not $B$" is: $  A \not \subseteq B$.
\\
With such interpretations, syllogistic may be treated as a calculus of equations [9], something that Boole seems to have had in mind. It reflects a rather impoverished \emph{Organon}. But if people  considered this sort of models as the true interpretation of Aristotle, there is no place for relations, and the opinion that he therefore is responsible for the lacunae of Euclidean axiomatics gets some support. 
\\
Scholars of Euclidean axiomatics such as de Risi [3], do not share this opinion: Aristotle, an alumnus of Plato's academy, where famously nobody entered without it, did know geometry. Prior Analytics contains a geometric proof of the equality of base angles in an equilateral triangle,(PA A 26), and  \textquotedblleft the principle of Aristotle" of Euclidean tradition is related to parallelity. 
\\
Posterior Analytics (Chapter 19), reflects Aristotle's understanding of \emph{recursion} as a mentally completed inductive definition of a concept. \textquotedblleft Mental completion" is hard to understand without a set-theoretic mindset, and it was a controversial issue for many commentators of Aristotle.\footnote{It was therefore only settled after the creation of set theory.}

\vspace {3 mm}

Logic definitely turned away from the Aristotelian tradition only at the turn to the 20th century. Bertrand Russell was an important mover in this. He had learned classical logic as a student but also had read up on Leibniz's attempts at reforming logic, critically. He asked himself several questions inspired by this reading, in particular one that is relevant here:
\\
\textquotedblleft  (1) Are all propositions reducible to the subject-predicate form?" [  8, p.13 ] \footnote{I'm grateful to Prof. V. de Risi for pointing me to this book in recent correspondence.}
\\
On the following pages of this book he proceeds to demonstrate by examples, that a logic adequate for mathematics cannot dispense with relations. Indeed later, in \emph{Principia Mathematica} they are a central ingredient.
\vspace{ 10 mm}

\newpage

\textbf{   3 \quad The \emph{Organon} in Combinatory logic}

\vspace {3 mm}

An adequate mathematical model for more of Aristotle's logic seems to be missing. This section describes my attempt to construct one.
\\
Logicians who have followed me to this place have long noticed that the \emph{ad hoc} formalism that I introduced in the discussion of Aristotle's family is in fact an extension of combinatory logic. Predicates, logical concepts and operations were added to the language by defining equations in these terms, combined with the combinatory application operation. This will now be turned into a formal calculus based on a language $\mathcal {E}$ which extends the language of combinatory logic.
\\
Before introducing this language, we choose a mathematical structure into which the language will be interpreted. Because it includes combinatory logic, we need a model for that.

\vspace {5 mm}

\textsc{ 3.1 The Modelling Structure }

\vspace { 3 mm}

The Language $\mathcal {E}$ of the modelling extends that of combinatory logic. 
\\
We first consider a fragment $\mathcal {E}$$_0$ of $\mathcal {E}$. It consists of expressions,  built from variables and constant predicates by the operations of application, the epsilon operator $\varepsilon_{x}$ and  the alpha operator $\alpha_{x}$ for variables $x$. These operations bind the variable $x$. We shall distinguish between objects that are called \emph{predicates} and predications.  A \emph{predication} may be obtained by applying a predicate $P$ to a variable $x$, \textquotedblleft predicating something about $x$", written $P \cdot x$.

\vspace{ 3 mm}

The basis for our interpretation of the language $\mathcal{E}$$_0$ is the graph-model of combinatory logic, Scott (1969), Plotkin (1972), Engeler (1981).
\\
Let $A$ be a non-empty set and define recursively
\\
$G_{0}(A) = A,  G_{n+1}(A) = A \cup \lbrace \alpha \rightarrow a: \alpha \ \text{finite or empty}, \alpha \subseteq G_{n}(A), a \in G(A) \rbrace$,
\\
where $\alpha \rightarrow a$ is a notation for the pair $\langle \alpha, a \rangle $. 
\\
The union of these $G_{n}(A)$ is denoted by $G(A)$. 
\\
The combinatory application operation is defined on subsets $M$ and $N$ of $G(A)$ as 
\\
$M \cdot N = \lbrace x :\alpha \rightarrow x \in M, \alpha  \subseteq N \rbrace $.

\vspace {3 mm}

With this interpretation of the application operation, the set $G(A)$ can be shown to be a model of combinatory logic. The elements of the model are the subsets of $G(A)$. We shall show below that the model satisfies the \emph{Comprehension Axiom of Combinatory Logic}: 
\\
For every expression $\varphi( x_1, \dots , x_n)$ built up by the application operation from the constants and variables (interpreted as subsets of $G(A)$),  there exists an element $M$ of the model such that
\\
$M \cdot x_{1} \cdot x_{2} \dots \cdot x_{n} = \varphi( x_1, \dots , x_n)$.
\\
The proof of this theorem actually produces an algorithm of comprehension to obtain $M$. Observe that all elements of $M$ have the form 
\\
$ \alpha_{1} \to (\alpha_{2} \to \dots ( \alpha_{n} \to a))$ with $\alpha_{i} \subseteq G(A)  \ \text {finite or empty, and } a \in G(A)$.
\vspace{ 3 mm}

For the next steps we shall rely on \emph{Sch\"{o}nfinkel's comprehension theorem}. He shows that if we have the \textquotedblleft combinators" $\boldsymbol {S}$ and $\boldsymbol {K}$ for which the  equations $\boldsymbol{S}PQR = PQ(PR)$ and $\boldsymbol {K}PQ = Q $ hold for all predicates $P, Q, R$, then there is the following conversion:
\\
 \quad  \emph{Comprehension Theorem of Combinatory Logic}
\\
For every combinatory expression $\varphi (x_1, \dots x_n)$ built up from constants and the variables $x_i$ there is a purely applicative expression $\psi(\boldsymbol {S}, \boldsymbol {K}) $ such that 
\\
$((((\psi(\boldsymbol {S}, \boldsymbol {K}) \cdot x_{1} )\cdot x_{2} )\cdots )x_{n} = \varphi(x_1, \dots, x_n)$.  
\\
For completing the proof we need only produce interpretations of the two constants and show, by inspection, that these conform to the equations, and thereby verify that by our interpretation we have in fact a model of combinatory logic. This is done in this author's 1981 paper on graph models, [5]. Here are the interpretations of the two combinators:
\\
$[\boldsymbol K] = \lbrace \lbrace a \rbrace \rightarrow (\emptyset \rightarrow a): a \in G(A)\rbrace $,
\\
$[\boldsymbol S] =  \lbrace ( \lbrace \tau \to ( \lbrace r_1, \dots,r_n\rbrace \to s)\rbrace \to ( \lbrace \sigma_{1} \to r_{1}, \dots \sigma_{n} \to r_{n} \rbrace \to
\\
(\sigma \to s))): n \geq 0,  r_1, \dots, r_n \in G(A), \tau  \cup \sigma_1 \cup \dots \cup \sigma_n = \sigma \subseteq G(A), \sigma  \quad \text {finite} \rbrace $.
 
 \vspace{ 5 mm}
 
\textsc {3.2 The Combinatory Interpretation of Categorical Predication}

\vspace { 3 mm}

\emph {Combinatory Predicates} are composed by the operation of application from predicate constants $C_j$ and variables $x_i$ to form expressions $ \varphi(C_1, \dots, C_m, x_1, \dots x_n)$. The constants $C_j$ are interpreted as subsets $[C_{j}]$ of $G(A)$, each variable $x_i$ ranges over a specific subset of $G(A)$. Their mention in $\varphi$ is usually suppressed.

\vspace { 3 mm}

Some predicates can be used for \emph{predications}: If the predicate $P$ is interpreted as $[P]$, a subset of $G(A)$, and $[P]$ is a set of elements of $G(A)$ of the form
\\
$ (\alpha_{1} \rightarrow (\alpha_{2} \rightarrow \dots (\alpha_{n} \rightarrow a))\ \text{with} \ a \in G(A),\alpha_{i} \subseteq G(A )\ \text{finite or empty}, i = 1, \dots, n $, 
\\
then $[P]$ can act as a predication $[P] \cdot[ x_{1}] \cdots [x_{n}]$ on these variables, interpreted as subsets $[x_i]$ of $G(A)$.
\\
\emph{Notation} : Where no ambiguity results we may omit the brackets on interpreted variables in the future.

\vspace { 3 mm}

The \emph{intuition behind this interpretation} of predication is that $[P]$ as a predication expresses some \emph{facts} about each subject-variable $x_i$. These facts are the extent to which $x_i$ conforms to the predicate $P$, the conformity being expressed by the corresponding sets $\alpha_i$. We call these facts \textquotedblleft \emph{attributes}".
\\
An interpretation of a predication is perhaps best illustrated by an example which we take from the family context of section 1. The interpretation of the parent  predicate $ [\textit {parent}(x, y)] $ is a set of expressions $(\alpha_1 \to ( \alpha_2 \to a))$ with $\alpha_1, \alpha_2 \subseteq G(A), a \in A$, and where $A$ is a set of people, each is present with the individual attributes. -- In distinction to section 1 we added the variables inside the the brackets for clarity, they relate the  variables $x, y$ to the sets $\alpha_1, \alpha_2$ in that order. Each set $\alpha_i$ is understood as a set of attributes: $\alpha_1$ of being a parent, $\alpha_2$ of being a child.
\\
The meaning of the predication $[\textit {parent}(x, y)]\cdot [x] \cdot [y]$ therefore is: \textquotedblleft $[y]$ is the set of people  in $A$ for whom $[x]$ is a set of parents". Specifically: $\alpha_1$ is the set of expressions $\{x\} \to \ ( \{z\} \to z)$ for x male, $\{y\} \to (\{z\} \to z)$ for female, and $\alpha_2$ consists of all $\{z\} \to z$ for the children z, with $x, y, z \in A$. The predication  produces the children of $x$ and $y$ if $x$ is male and $y$ female. 

\vspace {5 mm}

\emph { Categorical Predicates}, the analog to the categorical statements in the \emph{Organon}, arise from combinatory predicates by using \textquotedblleft \textit{for some}" and \textquotedblleft \textit{for all}", referring to the variables of a predicate. They constitute our language $\mathcal E$$_0$. We use $\varepsilon_x \varphi(x)$ to denote  \textquotedblleft some $x$ has $ \varphi(x)$" and $\alpha_x \varphi(x)$ to denote \textquotedblleft all $x$ have $ \varphi(x)$".
\\
Extending the modelling to the $\varepsilon$-operator is a bit subtle.  The term $\varepsilon_x \varphi(x)$ is to be interpreted as the result of a recursion in the sense of \textquotedblleft completed induction". Recall that the modelling of the language $\mathcal E$$_0$ is a process of finding denotations for elements of the language in a combinatory model. The modelling of $\varepsilon_x \varphi(x)$ involves the determination of an object $F$, a subset of $G(A)$, which has the property $\varphi(F)$. \textquotedblleft Recursion" means that such an object is already determined by an object $F_0$, the basis of recursion. This implies that the process of interpretation calls here for the choice of a particular object $F_0$,  which, as the case may be, is a challenge for the ingenuity of the modeller.
\\
 Given a unary predicate $P$, the object $[\varepsilon_{x} ( Px)]$ is therefore determined by the interpretation, which proposes an initial set $F_0 \subseteq G(A)$ with the property $F_0 \subseteq [P]\cdot F_0 $,  and yields
\\
 $[\varepsilon_{x} ( Px)]$ = $ \bigcup_{n} [P]^{n} F_0 = F$,  $[P]^{n}$ denoting the $n$-th iteration. 
\\ 
Then $F$ is a fixpoint of $[P]$, noting   $[P] \cdot F = F$.
\\
The finding of an appropriate $F_0$ is the cardinal point on which it turns whether or not the interpretation of the predication becomes vacuous, (see e.g. 4.1 below on the existence of a model for projective geometry). $F_0$ always exists, determined by the interpretation, in the worst case it is $F_0 = F = \emptyset $.
\\
The $\alpha$-operator is interpreted as 
\\
$[\alpha_x (Px) ]= [P] \cdot ext_{[P]}([x])$, where $ext_{[P]}([x])$ = $ \{a  :  \exists \alpha \to a \in [P] \}$ is the set of possible values for $[P]x$.
\\
The Aristotelian \textquotedblleft some $P$ are $Q$" thus translates into $[Q] \cdot [\varepsilon_{x}(Px)]$ and  \textquotedblleft all $P$ are $Q$" into $[Q] \cdot [\alpha_{x}(Px)] $.
\\
To extend these operations to $n$-ary predications we make another use of comprehension to separate out a specific variable in an expression $\varphi(x_1, \dots, x_n)$:
\\
$(\varphi_{j}( x_{1}, \dots, x_{n})\cdot x_1 \cdots x_{j-1} \cdot x_{j+1} \cdots x_{n}) \cdot x_{j} = \varphi( x_{1}, \dots, x_{n})\cdot x_{1} \cdots x_{n}$
\\
The $\varepsilon$ -operator and $\alpha$ -operators for $n$-ary predicates $[\varphi]$ are defined accordingly :
\\
$[\varepsilon_{x_j} ( \varphi_{j} (x_1, \dots, x_n))] \cdot x_{1} \cdots x_{j-1}  \cdot x_{j+1} \cdots x_{n}= F$, where, for $F_0$ given by the interpretation,
\\
$F =  \bigcup_{m}( [\varphi_{j} (x_1, \dots, x_n)]\cdot x_{1}\cdots x_{j-1} \cdot x_{j+1} \cdots x_{n})^{m} \cdot F_{0} $.
\\ê
$F$ is a set function with $n-1$ variables.
\\
$ [\alpha_{x_j} ([ \varphi_{j} (x_1, \dots, x_n))]\cdot x_{1} \cdots x_{j-1}  \cdot x_{j+1} \cdots x_{n} \\= [ \varphi_{j} (x_1, \dots, x_n)] \cdot x_{1}\cdots x_{j-1}  \cdot x_{j+1} \cdots x_{n}\cdot ext_{[\varphi_j]}([x_j])$, where 
\\
$ext_{[\varphi_{j}]} ([x_j])= \{ a : \exists(\alpha_{1} \to \dots \to (\alpha_{n} \to a)) \in [ \varphi_{j} (x_1, \dots, x_n))]\}$.

\vspace{3 mm}

\emph{Remark:} Two categorical statements \textquotedblleft no $P $ is $Q$" and \textquotedblleft some $P$ are not $Q$" are missing in $\mathcal{E}$$_0$. They are added in the next section in the context of negation. This is an expository choice. In fact they could have been added here separately, which would make $\mathcal{E}$$_0$  the full  categorical language.

 \vspace{ 5 mm}
 
\textsc{ 3.3 The Interpretation of Logical Connectives and Truth}

\vspace { 3 mm}

Predications as defined above are \textquotedblleft factual" interpretations, they produce a set of facts $[P] \cdot x_{1}  \cdots x_{n} $. Our interpretation of the language $\mathcal{E}$$_0$ resulted in a calculus of facts and as such cannot really be called a logical calculus. It lacks the logical connectives and judgements about the truth of a predication.
 \\
Predications lend themselves to logical composition by the connectives $\wedge, \vee$ and $\neg$. These constitute a language extension $\mathcal {E}$ of $\mathcal{E}$$_0$. The interpretation is extended to $\mathcal E$ recursively on the structure of the logical composition: the evaluation of
\\
$ [\varphi(x_1, \dots, x_n)] x_{1} \cdots x_{n} \wedge [\psi(x_1, \dots, x_n)] x_{1} \cdots x_{n}$ is \\ $[\varphi(x_1, \dots, x_n)] x_{1} \cdots x_{n} \cap [\psi(x_1, \dots, x_n)] x_{1} \cdots x_{n} $,
\\
correspondingly with $\vee$ and $\cup$, where we conformed the two predications to combined variables $x_1, \dots, x_n$ by comprehension. 
\\
For negation we set $ [\neg \varphi(x_1, \dots, x_n)] x_{1} \cdots x_{n} = ext_{[\varphi_{n}]}([x_n]) - [\varphi(x_1, \dots, x_n)] x_{1} \cdots x_{n}$.
\\
This concludes the definition of the logical predications. In particular, we can now express all the syllogistic statements \textquotedblleft some $P$ are $Q$", $\dots$, \textquotedblleft some $P$ are not $Q$", and our \emph{ad hoc} formalism in section 1 is thereby legitimised. 

\vspace {3 mm}

\emph{ Observation:} The operation of negation could have been added separately in the definition of the language $\mathcal {E}$$_0$ which would make it possible to add to it the two missing categorical statements \textquotedblleft no $P $ is $Q$" and \textquotedblleft some $P$ are not $Q$". The extended language $\mathcal{E}$ thus includes the full Aristotelian language of categorical statements.

\vspace {3 mm}

The intuition behind our \emph {Truth Definition} for a unary predicate $P$, modelled by a set of expressions $\alpha_i \to a$ is  that $[P]x$ is true if $x$ has all the attributes that are required by  $P$, that is $[P]x = \{ a : \alpha \to a \in [P] \}$. Correspondingly the truth definition for arbitrary predications in $\mathcal E$ is
\\
$[\varphi(x_1, \dots, x_n) x_{1} \cdots x_{n}$ is \emph{true}, denoted by $\top$, 
\\
if $[\varphi_{j}(x_1, \dots, x_n) x_{1} \cdots x_{j-1} \cdot x_{j+1} \cdots x_{n} \cdot x_{j} = ext_{[\varphi_{j}]}([x_j])$ for each $j$.
\\
$[\varphi(x_1, \dots, x_n) x_{1} \cdots x_{n}$ is \emph {false}, denoted by $\bot$, \\ if it is a proper non-empty subset of $ext_{[\varphi_{j}]}([x_j])$ for some $j$.
\\
In all other cases it is \emph {indeterminate}, denoted by $\triangle$.

\vspace { 5 mm}

\textsc{ 3.4 Relational Predications: From Factual to Logical Interpretation}

\vspace { 3 mm}

The \textquotedblleft facts" produced by predications $[P]\cdot x_{i} \cdots x_{n}$ in $\mathcal{E}$ do not reflect the actual relations between the arguments $x_i$. This can be accomplished by introducing \emph{relational predications} $[P_R]$:
If $[P]$ consists of elements$ (\alpha_{1} \to \dots \to (\alpha_n \to a))$ with $a \in G(A), \alpha_i \subseteq G(A) \ \text{finite and nonempty}$, then
\\
$[P]_R\cdot x_{1} \cdots x_{n}$ is the set of all $\langle a_1, \dots, a_n \rangle: \exists (\alpha_{1} \to \dots \to (\alpha_n \to a)) \in [P] $
\\ 
 such that $ a_i \subseteq x_i \subseteq \alpha_i, i = 1, \dots, n $. The subscript $R$ is tacitly understood in the following.

\vspace{ 3 mm}

Let $\mathcal {E}$ be extended to  $\mathcal{E}$$^{\Lambda}(C_1, \dots, C_m)$ by introducing a valuation of the relational predicate constants $C_1, \dots, C_m$ of $\mathcal{E}$.

\vspace {3 mm}

The \emph{factual interpretation} of relational predications distinguishes items that verify or falsify it. This distinction is based on a valuation $\Lambda_i$ for each one
of the constant predicates $[C_i]$. We denote the valuation of a tuple $\langle a_1, \dots, a_n \rangle$ by $\langle a_1, \dots, a_n \rangle ^{\Lambda_i} = \langle a_1, \dots, a_n \rangle ^{\top}$ if it is a \emph{verifying} fact, by $\langle a_1, \dots, a_n \rangle ^{\bot}$ if it is \emph{falsifying}.
\\
For a $n$-ary predicate constants $[C_i]$ the facts that are valued by the valuation $\Lambda_i$ are the tuples  $\langle a_1, \dots, a_n \rangle \in D_i$ , where
\\
  $D_i $ = $\lbrace \langle a_1, \dots, a_n \rangle : (\alpha_{1} \to \dots (\alpha_{n-1} \to  a_j)) \in [(C_i)_j], j = 1, \dots, n \rbrace$ 
\\
We define the corresponding predication  $[C_i]^{\Lambda_i}$ by  the set of objects 
\\
$ ( \alpha_{1} \to \dots (\alpha_{n} \to \langle a_1, \dots, a_n \rangle         )^{\Lambda_{i}}$ where
\\
$ (\alpha_{1} \to (\alpha_{2} \to \dots  \to (\alpha_{n} \to a_{n}))) \in [C_i], \langle a_1, \dots, a_n \rangle \in D_i $.
\\
As a result, the factual interpretation of the relational predication $C_i$ is the set function $[C_i]^{\Lambda_i} \cdot x_{1} \cdots x_{n}$ which produces a set containing elements $ \langle a_1, \dots, a_n \rangle ^{\top}, \\ \langle a_1, \dots, a_n \rangle ^{\bot}$ and $ \triangle$, (for the cases where the predication returns the empty set on the given inputs $x_1, \dots, x_n$). 

\vspace {3 mm}

The factual interpretation of the language $\mathcal{E}$$^{\Lambda}(C_1, \dots, C_m)$ is based on the valuation $\Lambda_{i}$ for each $C_i$, and then extended over logical composition, the $\varepsilon$- and $\alpha$-operations and the quantifiers $\exists$ and $\forall$ as follows:   The valuation $\Lambda$ maps each predication into a set of the above objects. We represent each such set as a propositional formula consisting of the valued tuples $\langle a_1, \dots, a_n \rangle$. The valuation $\Lambda$ assigns to each of them the truth-value \emph{true} or \emph{false} as indicated by the superscripts. To obtain the truth-value of a constant predication, the set produced by it is interpreted as the conjunction of these elements as a formula in a propositional logic. The presence of $\triangle$ in a propositional formula assigns to it the value \textit{indeterminate}.
\\
This interpretation is then extended as follows to predications obtained by the logical operations:
\\
The conjunction and disjunction of predications $[\varphi(x_1, \dots x_n)]^{\Lambda}\cdot x_{1} \cdots x_{n}$ and \\ $[\psi(y_1, \dots y_m)]^{\Lambda}\cdot y_{1} \cdots y_{m}$
are conjunctions respectively disjunctions of the corresponding propositional expressions. The interpretation of negation is obtained by  inverting all $\top$ to $\bot$ and all $\bot$ to $\top$ in the valuations of the tuples $\langle a_1, \dots, a_n \rangle ^{\top},\\ \langle a_1, \dots, a_n \rangle ^{\bot}$. 

\vspace { 3 mm}

The factual interpretation of the $\varepsilon$- and $\alpha$ operators on a predicate $\varphi(x_1, \dots, x_n)$ is obtained as
\\
$[\varepsilon_{x_j} (\varphi(x_1, \dots, x_n)]^{\Lambda}\cdot x{_1} \cdots x_{j-1} \cdot x_{j+1} \cdots x_{n}  \\
= \lbrace \langle a_1, \dots, a_{j-1},a_{j+1},\dots, a_n,a_j \rangle^{\top}: \\a_j \in  [\varepsilon_{x_j} (\varphi_{j}(x_1, \dots, x_n)]^{\Lambda}\cdot x{_1} \cdots x_{j-1} \cdot x_{j+1} \cdots x_{n} \cdot x_{j}\rbrace $,
\\
$[\alpha_{x_j} (\varphi(x_1, \dots, x_n)]^{\Lambda}\cdot x_{1} \cdots x_{j-1} \cdot x_{j+1} \cdots x_{n} \\ 
= \{ \langle a_1, \dots, a_{j-1},a_{j+1},\dots, a_n,a_j\rangle^{\top} : \\ a_j \in [\varphi_{j}(x_1, \dots, x_n)]^{\Lambda}\cdot x_{1} \cdots x_{j-1}\cdot x_{j+1} \cdots x_{n}\cdot ext_{[\varphi_j]}([x_{j}]) \} 
\\
\cup \{ \langle a_1, \dots, a_{j-1},a_{j+1},\dots, a_n,a_j\rangle^{\bot} : \\ a_j \in [ \neg\varphi_{j}(x_1, \dots, x_n)]^{\Lambda}\cdot x_{1} \cdots x_{j-1}\cdot x_{j+1} \cdots x_{n}\cdot ext_{[\varphi_j]}([x_{j}]) \}$.

\vspace{ 3 mm}

\emph {Existential and universal quantifiers} can now be introduced as follows:
\\
$[\exists x_{j} \varphi(x_1, \dots, x_n)]^{\Lambda} \cdot x_1 \cdots x_{j-1} \cdot x_{j+1} \cdots x_{n}\\
= [\varepsilon_{x_j} (\varphi(x_1, \dots, x_n)]^{\Lambda}\cdot x{_1} \cdots x_{j-1} \cdot x_{j+1} \cdots x_{n} $,
 \\
 $[\forall x_{j} \varphi_{j} (x_1, \dots, x_n)] ^{\Lambda} \cdot x_1 \cdots x_{j-1} \cdot x_{j+1} \cdots x_{n}\\
=  [\alpha_{x_j} (\varphi(x_1, \dots, x_n)]^{\Lambda}\cdot x_{1} \cdots x_{j-1} \cdot x_{j+1} \cdots x_{n}$.

 \vspace { 3 mm}
 
 Altogether, we have now extended the factual interpretation of predications to all of the language $\mathcal{E}$$^{\Lambda}(C_1, \dots, C_m)$. Each predication produces a propositional formula containing elements of the form $\langle a_1, \dots, a_n \rangle ^{\bot} ,\langle a_1, \dots, a_n \rangle ^{\top}$ and $\triangle$. This completes the factual interpretation of  $\mathcal{E}$$^{\Lambda}(C_1, \dots, C_m)$; it assigns a propositional formula of verifying, falsifying tuples and $\triangle$-s. 
 \\
 The \emph{truth-value} of the predication is obtained by evaluating the tuples of the propositional components as true, false and indeterminate as above.

\newpage

\textbf{   4 \quad Discussion}

\vspace {3 mm}

\textsc{4.1 Projective Geometry}

\vspace { 3 mm}

Looking for an example to discuss our interpretation $\mathcal {E}$ and $\mathcal{E}$$^{\Lambda}$ of the \emph{Organon}, recall Plato's advice in the \emph{Republic}, (Chapter VII), that with geometry you can educate the mind.
\\
We take \textquotedblleft geometry" here as a first-order theory  in mathematical logic, and for simplicity restrict to projective geometry.

\vspace { 3 mm}

The first-order models in mathematical logic are relational structures \\  $ \langle A, R_1, \dots, R_n, c_1, \dots,c_n \rangle$ with relations $R_i \subseteq A^{k_i}$ and individual constants $c_j \in A$. In our modelling the relations $R_i$ correspond to combinatory constants that denote $k_i$-ary predicates; the constants $c_j$ also correspond to combinatory constants. These objects are then interpreted as subsets of $G(A)$ according to the above definitions. 
\\
As an example, consider a mathematical structure such as a projective plane, understood as a set $P$ of \textquotedblleft points", $L$ as a set of \textquotedblleft lines" with the binary relation of \textquotedblleft incidence" $\textit{Inc}\subseteq P \times L$, where $P \cap L = \emptyset$. These are subject to some axioms such as: \textquotedblleft For any two points there is a unique line on which they are incident."
\\
The combinatory model is based on this $A$ for the construction of $G(A)$. In $\mathcal {E}$ the incidence relation is interpreted as 
\\
$[\textit{Inc}]^{\Lambda} = \{ (\{ p \} \to (\{ p\} \to l)) \to (\{ p\} \to l): p \in P, l \in L, \langle p, l \rangle \in Inc \}$.
\\
In $\mathcal{E}$$^{\Lambda}(Inc)$ it would be
\\
$[\textit{Inc}]^{\Lambda} = \{ (\{ p \} \to (\{ p\} \to l)) \to (\{ p\} \to \langle p, l \rangle): p \in P, l \in L, \langle p, l \rangle \in Inc \}$
\\
An equality predicate is needed here only for points and lines and can therefore be viewed as a binary predicate constant with the interpretation
\\
$[eq] = \{ (\{ x \} \to (\{ y\} \to y)) \to ( \{ y \} \to y): x= y, x, y \in P \cup L \} $ in $\mathcal {E}$, 
\\
with the corresponding the relational predication
\\
$[eq]^{\Lambda} = \{(\{ x \} \to (\{ y\} \to y)) \to (\{ x\} \to \langle x, y \rangle): x= y, x, y \in P \cup L \} $ in $\mathcal{E}$$^{\Lambda}$ .

\vspace {3 mm}

Using the quantifiers introduced earlier, the above axiom is interpreted as
\\
$\forall x_{1} \forall x_{2} ( [eq]^{\Lambda}x_{1} x_{2}  \vee \exists y  ( [\textit{Inc}]^{\Lambda}x_{1}y \wedge [\textit{Inc}]^{\Lambda}x_{2}y \wedge \forall z ( [\neg \textit{Inc}]^{\Lambda}x_{1}z \vee [\neg \textit{Inc}]^{\Lambda}x_{2}z \vee [eq]^{\Lambda}yz )))$.
\\
with the parameter $[\textit{Inc}]$. -- The other axioms would be represented in the same fashion and combined into a logical predicate denoted by $\pi(\textit{Inc})$. It then turns into an exercise of inventiveness in finding the required fixpoints, and of formal persistence to verify that it returns the value \emph{true} on some given model of projective plane geometry, that is for a given binary relation of incidence, e.g. the \emph{Fano plane} of seven points and seven lines. 

\vspace { 3mm}

Of course, projective geometry itself may be considered as a defined predicate, obtained by using a suitably defined recursion on a predication: Let the variable $X$ be substituted for $\textit{Inc}$ in the axiom-expression $\pi (\textit{Inc})$ . The recursion equation $ X = \varepsilon_{X}\pi(X))$ defines the geometric concept of a projective plane. The modelling of this formula starts with fixing on a given set of points and lines, these may come from defniitions in terms of finite fields or $\mathbb {Q}, \mathbb{R}$ or $\mathbb{C}$. Here again the modelling consists in the finding of a fixpoint, which is the crucial matter in obtaining an actual model.  If the language has more than one predicate constant, e.g. one for betweenness, one would of course use joint recursion for incidence and betweenness.
\\
\emph{ Remark }: Perhaps it is worth mentioning that this use of the predicate $\pi(X)$ is an example of introducing additional predicate constants with which one may  extend the language to encompass additional concepts. This corresponds to the familiar way to define new mathematical concepts and structures.

\vspace{ 10 mm}

\textsc {4.2 Facts and Thoughts }

\vspace { 3 mm}

The restriction of \textquotedblleft facts" to the the basic set $A$ as the extension of the variables is natural in simple contexts  like the Aristotle family. We also used it in the context of projective geometry treating of "facts" about points and lines. But this is too restrictive even in this case: Consider the notion of point-functions, for example projectivities. The object $f$ is a point-function if $fpq_1 = fpq_2$ for all $q_1,q_2 \in P$. This can be expressed by a predication $[\textit{fun}]$ defined by
\\
$[\textit{fun}]f = [eq] \cdot \alpha_{q_1}(fpq_1) \cdot \alpha_{q_2} (fpq_{2})$,
\\
using the predication $[eq]$ from above. Observe that $[\textit{fun}]$ is an element of $G_2(A)$.

\vspace {3 mm}

The $\varepsilon$-operator also creates \textquotedblleft facts": Consider the line connecting two point $p_1, p_2$, 
\\
expressed by $\varepsilon_{l} ([eq]\cdot [\textit{Inc}](p_1l) \cdot (p_2l))$,
\\
which is a function on $p_1, p_2$, called a \emph{Skolem-function} in logic. It is an element of $G_2(A)$.
\\
A \textquotedblleft Skolem-function" $f$ is thus the result of a recursion, a concept that we traced back to Aristotle's notion of a mentally completed definition of induction. It is therefore legitimate to call this object a \emph{thought}. -- Anyway, I would have preferred \emph{\textquotedblleft thoughts"} over \textquotedblleft \emph{predicates}" and \textquotedblleft \emph{predications}". The latter are naturals in the Aristotelean context. But I see predicates as thoughts, as sets or patterns of small and big notions: The predication $[P] \cdot x$ is perceived as applying a thought $P$ to a thought $x$, checking to what extent the thought $[P]$ applies to $x$. This perception is the background of my modelling and is connected to my work on neural algebra which treats of thoughts as patterns of firing neurons.

\vspace{ 5 mm}

\textsc {4.3 Algorithmic Logic }

\vspace { 3 mm}

Let me not forget my own brain-child, algorithmic logic, [4], now approaching retirement age after a long career. It treats of algorithmic properties of structures. Its predicates are of the form $\pi(x_1, \dots, x_n)$ which denotes a program $\pi$ with the input variables $x_1, \dots x_n$ of elements of the structure. Since combinatory logic can deal with the notions of computation and termination, there is an important point of contact here which merits elaboration.
\\
Programs $\pi(x_1, \dots, x_n)$ are composed of individual instructions, namely assignments of the form $z := f(x, y)$, decisions such as $ x < y$. These correspond, loosely speaking, to our predicate constants. Program statements are composed by successive execution  $(\pi_1(x_1, \dots, x_n)) ; $  $(\pi(_2x_1, \dots, x_n))$ and recursion. These essentially correspond to composition and the $\varepsilon$-operation. Finally, the factual interpretation $[ \pi(x_1, \dots, x_n)]^{\Lambda} \cdot{x_1} \cdots x_{n}$ of a program is the so-called denotational semantics of the program. The valuation $\Lambda$ is understood as the valuation of tuples of elements of the relations and functions in some relational structure. The program statement $[\pi(x_1, \dots, x_n)]^{\Lambda} x_{1} \cdots x_{n}$ evaluates to the result of executing the program on the input assignment.
\\
For the logic of programs, the \textquotedblleft \emph{algorithmic logic}", we chose to evaluate a program as \textquotedblleft true" in a relational structure if it halts on all inputs.
\vspace {10 mm}

\textbf { 5 \quad  Conclusions and Scholia}

\vspace{ 3 mm}

\textsc { 5.1 On the discussion of Aristotle's Relations}

\vspace{ 3 mm}

This discussion may be conducted in the interpretation of the language  $\mathcal{E}$$^{\Lambda}$ in section 3.4 above which captures, extends and completes the \emph{ad hoc }formalism used in our motivational section 1. The example fits nicely into this framework: The set $A$ lists the names of the members of the family. The binary relations  for motherhood, fatherhood and marriage are lists of pairs of names, and singletons for being male or female. Thus, the interpretation $[father]^{\Lambda} $ for fatherhood is a case of binary relations which are represented in the form
\\
$\{\{ x \} \to (\{ y\} \to y)) \to (\{ x\} \to \langle x, y \rangle)$.
\\ 
Therefore fatherhood in the Aristotle family would be represented by  a set containing the substitution instances 
\\
$ x :=$ \textit{Niarchus1}, $y:=$ \textit{Aristotle} and $x:=$ \textit{Niarchus1}, $y:=$ \textit{Arimnestus} 
\\
as well as $x:=$ \textit{Aristotle}, $y :=$ \textit{Niarchus2}.
\\
A more adequate modelling of fatherhood would need more attributes from of the vital statistics of the family members, e.g. profession, place and date of birth, etc. to exclude false claims of fatherhood. Above, we have artificially distinguished the two people called \textquotedblleft Niarchus" by adding the distinction to the names.
\\
In this setting it becomes clear that Aristotle had the conceptual means to conceive of, and formally treat, relations. But the primary goal of the \emph{Organon} did not require this.

\vspace {3 mm}

This concludes our search for Aristotle's relations. The missing relations of Euclid, a blemish on his axioms, were not noticed at the time because geometry was understood by Aristotle as describing geometric properties of lines and other objects unquestioned as continuous. As pointed out by de Risi in [2 ] this was only put into question in the 16th century and not fully received into the understanding of space till the 19th century.\footnote{The de Risi references were pointed pout to me by Prof.M.Beeson.}

\vspace { 3 mm }

\textsc { 5.2 Scholia}

\vspace {3 mm}

What have we learned from our combinatory experiment on the \emph{Organon} ?
\\
I:  The objects of syllogisms, the categorical statements, were interpreted as elements of the language $\mathcal {E}$ which expands the set of terms of combinatory logic. This combinatory interpretation thus becomes a model of syllogistics comparable to the models of \L ucasiewicz and others.
\\
II: The categorical statements themselves are statements about properties of facts, traditionally about individual facts. This aspect is captured by our language $\mathcal {E}$.
\\
III: The logical interpretation of the language $\mathcal{E}$$^{\Lambda}(C_1, \dots, C_m)$ is based on facts about relations as $n$-ary predications. This is probably  the most adequate rendering of my understanding of an extension of the \emph{Organon} if it were to include relations. 
\\
IV: Logical interpretations of predications turns them into \emph{Judgements}. Our distinction between factual and logical interpretation recalls the famous distinction between judgements \emph{de re} and judgements \emph{de dicto} going back to 12th century Scholastics, when Peter Abelard derived \emph{de dicto} from \emph{de re}, as we do.
\footnote{This was strongly contested at the time by Bernard de Clairvaux, who maintained that the dogmata of the Church (to which Abelard addressed himself) were \emph{de dicto} statements and not to be made dependent on judgements \emph{de re}. Bernard was sainted, Abelard not.}
\\
V: Our modelling addressed the \emph{semantics} of $\mathcal{E}$ and $\mathcal{E}$$^{\Lambda}$ and not its \emph{deductions}. But this is another chapter.

\vspace{10 mm}

\textbf{Apologia and Dedication} 

\vspace {3 mm}

Who knows how Aristotle would react to my experiment. I picture  him and Euclid as  little Raffael-angels looking down on our travails with mischievous interest.

\vspace { 3 mm}

 My interpretation is based on my individual reading of his \emph{Organon}. Individual does not mean indivisible; I may have had two minds about a number of things. Also, while I avoided the pitfalls of turning the logical connectives into predicates and moreover turn syllogisms themselves into statements, I do not propose to justify this here.

\vspace {3 mm}

Finally, I claim the privilege of a very old man and refrain from the labours of performing all the verifications and of following up my own suggestions. Let the friendly reader smile and forgive. Here, thanks are due to Michael Beeson whose detailed comments on the manuscript were very helpful.

\vspace { 3 mm}

I dedicate this little essay to all my friends and students who troubled themselves to remember me and my birthday at the symposium in Zurich, early 2020. And to all those that could not be present because of the distances that destiny put between us, in particular to my late friend and colleague Ernst Specker whom this symposium was meant to honour too. Special mention goes to the three organisers Gerhard Jaeger, Reinhard Kahle and Giovanni Sommaruga,  the sponsors that they found, and to my \emph{alma mater}, the ETH, for its hospitality.

\vspace {15 mm}

 \textsc {References}
\\
 {[1] John Corcoran. Aristotle's Prior Analytics and Boole's Laws of Thoughts. \emph{History and Philosophy of Logic}, 24: 261-288,2003
 \\
 {[2] Vincenzo de Risi. Francesco Patrizi and the New Geometry of Space}. In: Koen Vermeir and Jonathan Regier, editors. \textit{Boundaries, Extents and Circulations: Space and Spaciality in Early Modern Natural Philosophy}, number 41  in Studies in History and Philosophy of Science, chapter 3, pages 55 - 100. Springer International Publishing, Switzerland, 2016.
  \\
 {[3]) Vincenzo de Risi. The Development of Euclidean Axiomatics. \textit{ Archiv for the History of the Exact Sciences}, 70: 561-676, (2016)..
 \\
 {[4]) Erwin Engeler. Algorithmic Properties of Structures. \textit{Mathematical Systems Theory}, 1: 183-193, 1967.
 \\
 {[5]) Erwin Engeler. Algebras and Combinators. \textit{ Algebra Universalis}, 13:  389-392, 1981.
 \\
{[6]) Jan \L ukasiewicz. \textit{Aristotle's Syllogistics from the Standpoint of Modern Logic}. Clarendon Press, Oxford, second edition,1957.
 \\
{[7]) Octavius Freire Owen. \textit{The Organon, or Logical Treatisies of Aristotle with an Introduction to Porphyry}, (two volumes). Henry G. Bohn, London, 1853.
\\
{[8]) Bertrand Russell. \textit{A Critical Exposition of the Philosophy of Leibniz}. Cambridge Univ.Press, Cambridge, 1900.
 \\
{[6])  John C. Shepherdson. On the Interpretation of Aristotelean Syllogistics. \textit{Journal of Symbolic Logic}, 21: 137-147, 1956
 \\
{[10]) Robin Smith. \textit{Prior Analytics}. Hackett Publ.Co., Indianapolis Cambridge, 1989.

\nocite{*}
\bibliography{engeler}
\bibliographystyle{plain-annote}

\end{document}